\documentclass[12pt]{amsart}

\def \fg{{\mathfrak g}}
\def \fh{{\mathfrak h}}
\def \fk{{\mathfrak k}}

\def \fa{{\mathfrak a}}
\def \fb{{\mathfrak b}}
\def \fl{{\mathfrak l}}

\def \fn{{\mathfrak n}}
\def \ft{{\mathfrak t}}

\def \sl{{\mathfrak s}{\mathfrak l}}
\def \su{{\mathfrak s}{\mathfrak u}}

\def \fm{{\mathfrak m}}
\newcommand{\beqa}{\begin{eqnarray*}}
\newcommand{\eeqa}{\end{eqnarray*}}

\def \la{{\langle}}
\def \ra{{\rangle}}
\def \R{{\mathbb R}}
\def \P{{\mathbb P}}
\def \C{{\mathbb C}}

\def \i{{\rm i}}
\def \fd{{\bullet}}
\def \proof{{\noindent {\it Proof.} \ \ }}

\def \dim{{\rm dim}}

\newtheorem{theorem}{Theorem}[section]
\newtheorem{lemma}[theorem]{Lemma}
\newtheorem{proposition}[theorem]{Proposition}

\newtheorem{remark}[theorem]{Remark}
\newtheorem{corollary}[theorem]{Corollary}

\newtheorem{nota}[theorem]{Notation}

\numberwithin{equation}{section}

\begin{document}

\baselineskip=16pt

\title[Poisson structures on flag manifolds]{Poisson structures on complex flag manifolds
associated with real forms}

\author[P. Foth and J.-H. Lu]{Philip Foth and Jiang-Hua Lu}

\address{Department of Mathematics, University of Arizona, Tucson, AZ 85721-0089
\newline Department of Mathematics, The University of Hong Kong, Pokfulam, Hong Kong}

\email{foth@math.arizona.edu, jhlu@maths.hku.hk}

\subjclass{53D17}

\keywords{Lie groups, real forms, flag varieties, Poisson structures, symplectic leaves.}

\date{July 17, 2003}

\begin{abstract}
For a complex semisimple Lie group $G$ and a 
real form $G_0$ we define a Poisson structure on the variety of Borel subgroups 
of $G$ with the property that all $G_0$-orbits in $X$ as well as all  
Bruhat cells (for a suitable choice of a Borel subgroup of $G$) are Poisson 
submanifolds. In particular, we show that every non-empty intersection of a 
$G_0$-orbit and a Bruhat cell is a regular Poisson manifold and we
compute the dimension of its symplectic leaves. 
\end{abstract}

\maketitle

\baselineskip=12pt
\begin{flushright}{\it Dedicated to Alan Weinstein on \\
the occasion of his 60th Birthday.}
\end{flushright}
\baselineskip=16pt

\section{Introduction.}

Let $G$ be a connected and simply connected complex semisimple Lie group with Lie
algebra $\fg$, and let $X$ be the variety of Borel subalgebras of $\fg$. In this paper we 
use a real form $\fg_0$ of $\fg$ to define a Poisson structure on $X$. This Poisson structure
depends on a choice of a Borel subalgebra $\fb$ of $\fg$ such that $\fg_0 \cap \fb$ is a 
maximally compact Cartan subalgebra of $\fg_0$. Instead of dealing with each real form 
individually, we fix a Borel subalgebra $\fb$ of $\fg$ and a Cartan subalgebra 
$\fh \subset \fb$. Then, as is shown in \cite{Knapp}, a real from $\fg_v$
of $\fg$ can be constructed from each Vogan diagram $v$ for $\fg$ such that
 $\fg_v \cap \fb$ is a maximally compact Cartan subalgebra of $\fg_v$. The 
corresponding Poisson structure on $X$ is denoted by $\Pi_v$. 

Let $G_v$ be the real form of $G$ corresponding to $\fg_v$, and let $B$ be the 
Borel subgroup of $G$ with Lie algebra $\fb$. The Poisson
structure $\Pi_v$ has the property that each $G_v$-orbit as well as each $B$-orbit in $X$ 
is a Poisson submanifold. The $B$-orbits in $X$ will be referred to as the Bruhat cells. 
We compute the rank of $\Pi_v$. In particular, if $G_v$-orbit 
${{\mathcal O}}$ meets a Bruhat cell ${{\mathcal C}}$, they intersect transversally, 
and we find that all the symplectic leaves in ${{\mathcal O}} \cap {{\mathcal C}}$ have 
the same dimension, so ${{\mathcal O}} \cap {{\mathcal C}}$ is a regular Poisson manifold. 
Moreover, we show that all symplectic leaves in each connected component of 
${{\mathcal O}} \cap {{\mathcal C}}$ are translates of each other by elements of
a Cartan subgroup of $G_v$. We also show that the $G_v$-invariant Poisson cohomology for
each open $G_v$-orbit in $X$ is isomorphic to the de Rham cohomology of $X$. 

Results similar to those presented here for the full flag 
manifold $X=G/B$ are also valid for a partial flag manifold $G/P$, where $P$ is a 
parabolic subgroup of $G$ containing $B$. We will treat these more general cases as well 
as some further properties of $\Pi_v$ in a future paper. 

Throughout this paper, if $V$ is a set  and $\sigma$ is an involution on $V$, we will 
use $V^\sigma$ to denote the fixed point set of $\sigma$ in $V$.

\section{Real forms of $\fg$ and Vogan diagrams}
\label{sec_gv}

Let $\fg$ be a complex simple Lie algebra. In this section we recall the classification of
real forms of $\fg$ by Vogan diagrams. Details can be found in \cite[Chapter 6]{Knapp}.

Suppose that $\fg_0$ is a real form of $\fg$ and that $\tau_0$ is the corresponding 
complex-conjugate linear involution on $\fg$. Let $\theta_0$ be a Cartan involution of $\fg_0$,
and let $\fh_0$ be a $\theta_0$-stable maximally compact Cartan subalgebra of $\fg_0$. 
Set $\ft_0 = \fh_{0}^{\theta_0}$ and $\fa_0 = 
\fh_{0}^{-\theta_0}$ so that $\fh_0 = \ft_0 + \fa_0$. Let $\gamma_0$ be the 
complexification of $\theta_0$. Then the Cartan subalgebra $\fh = \fh_0 + \i \fh_0$
of $\fg$ is $\gamma_0$-stable. Let $\Delta$ be the root system for $(\fg, \fh)$. Since 
$\fh_0$ is a maximally compact  Cartan subalgebra of $\fg_0$, there exists $x_0  \in \i\ft_0$
that is regular for $\Delta$. Define the subset $\Delta^+$ of  positive roots  in $\Delta$ by
$\alpha \in \Delta^+$ if and only if $\alpha(x_0) > 0$. Then $\gamma_0 (\Delta^+) = \Delta^+$.
Let $\Sigma \subset \Delta^+$ be the set of simple roots in $\Delta^+$. Then $\gamma_0(\Sigma)
=\Sigma$, so $\gamma_0$ gives rise to an involutive automorphism of the Dynkin diagram of $\fg$.
 Let ${{\mathcal I}}$ be the set of non-compact imaginary simple roots.
The Vogan diagram of $\fg_0$ associated to the triple
$(\theta_0, \fh_0, \Delta^+)$ is the Dynkin diagram $D(\fg)$ of $\fg$
together with an involutive automorphism $\gamma_0$ on $D(\fg)$ and the 
vertices corresponding to the simple roots in ${{\mathcal I}}$ painted black. 

In general, a  Vogan diagram for $\fg$ is defined to be a triple $(D(\fg), d, {{\mathcal I}})$,
where $D(\fg)$ is the Dynkin diagram  of $\fg$, $d$ is an involutive automorphism of $D(\fg)$, 
and  ${{\mathcal I}}$ is a subset of vertices of
$D(\fg)$ such that $d(\alpha) = \alpha$ for each $\alpha \in {{\mathcal I}}$.
Every Vogan diagram for $\fg$ comes from a real form of $\fg$ (see below), although 
two different Vogan diagrams can come from isomorphic real forms. 
A non-redundant list of Vogan diagrams with the corresponding
isomorphism class of real forms for all simple Lie algebras is given in \cite{Knapp}.
Every Vogan diagram in the list in \cite{Knapp} is {\it normalized} in the sense
that at most one vertex is painted black. 

For the purpose of defining Poisson structures on the variety of Borel subalgebras
of $\fg$, we now recall the explicit construction of a real form of $\fg$ from a
Vogan diagram \cite[Theorem 6.88]{Knapp}.
We need to fix the following data for $\fg$.

Choose a Cartan subalgebra $\fh$ of $\fg$ and let $\Delta$ be the 
root system for $(\fg, \fh)$. 
Fix a choice of positive roots $\Delta^+$ and let $\Sigma$ be the basis of simple roots. 
Let $\ll, \, \gg$ be the Killing form of $\fg$ and let root vectors $\{E_\alpha:
\alpha \in \Delta\}$ be
chosen such that $[E_\alpha, E_{-\alpha}] = H_\alpha$ for each $\alpha \in \Delta^+$,
where $H_\alpha$ is the unique element of $\fh$ defined by $\ll H,
H_\alpha\gg=\alpha(H)$ for all $H\in\fh$, and such that the numbers $m_{\alpha, \beta}$
given by $[E_\alpha, E_\beta] =  m_{\alpha, \beta}E_{\alpha+\beta}$ when
$\alpha+\beta\in \Delta$ are real. Define a compact real form $\fk$ of $\fg$ 
as 
\[
\fk = {\rm span}_{{\mathbb R}} \{\i H_\alpha, \, X_\alpha:=E_\alpha-E_{-\alpha},
\, Y_\alpha:=\i(E_{\alpha}+E_{-\alpha}) \}\ ,
\]
and let $\theta$ be the complex conjugation of $\fg$
defining $\fk$. If $d$ is an involutive automorphism of the Dynkin diagram of $\fg$, define
$\gamma_d$ to be the unique automorphism of $\fg$ satisfying
$\gamma_d(H_\alpha)=H_{d(\alpha)}$ and $\gamma_d(E_{\alpha})=E_{d(\alpha)}$ 
for each simple root $\alpha$. 

Given a Vogan diagram $v$ for $\fg$, not necessarily normalized, with the involutive diagram
automorphism $d$, let $t_v$ be the unique element in the adjoint group of $\fg$ such that
$$
{\rm Ad}_{t_v}(E_\alpha)=\begin{cases}
E_\alpha \mbox{  if } \alpha \mbox{ is a blank vertex in $v$} \\
- E_\alpha \mbox{ if } \alpha \mbox{ is a painted vertex in $v$} 
\end{cases} 
$$ 
Define a complex conjugate linear involution 
$$
\tau_v:={\rm Ad}_{t_v}\circ \gamma_d\circ\theta.
$$ 

\begin{nota}
\label{nota_v}
{\rm
We use $\fg_v = \fg^{\tau_v}$ to denote the
real form of $\fg$ defined by $\tau_v$. 
Set $\theta_v = \theta|_{{\mathfrak g}_v}$. Then $\theta_v$ is a Cartan 
involution of $\fg_v$,
and $\fh^{\tau_v}$ 
is a $\theta_v$-stable maximally compact Cartan subalgebra of $\fg_v$, with
$\fh = \fh^{\tau_v}  + \i\fh^{\tau_v}$. The complexification of $\tau_v$ is 
\begin{equation}
\label{eq_gamma_v}
\gamma_v := \tau_v \theta = \theta \tau_v = {\rm Ad}_{t_v} \gamma_d.
\end{equation}
Since $\gamma_v(\Delta^+) = \Delta^+$,
the Vogan diagram of $\fg_v$ associated to the triple $(\theta_v, \fh^{\tau_v}, \Delta^+)$
is $v$.
}
\end{nota}

One of the advantages of introducing the real form $\fg_v$ is as follows. We say that
a real subalgebra $\fl$ of $\fg$ is {\it Lagrangian} if its real dimension is equal 
to the complex dimension of $\fg$ and if ${\rm Im} \ll x_1,x_2\gg =0$ for all $x_1, x_2\in\fl$.
A decomposition $\fg=\fl_1+\fl_2$ is called a 
{\it Lagrangian splitting} if both $\fl_1$ and $\fl_2$ are Lagrangian. Let
$\fn$ be the subalgebra of $\fg$ spanned by the set of all positive root vectors for $\Delta^+$.
The following fact is easy to prove. 

\begin{lemma}
\label{lem_splitting}
 Let $\fl_d:=\fh^{-\tau_v}+\fn$. Then
$\fg=\fg_v+\fl_d$ is a Lagrangian splitting of $\fg$.  
\end{lemma}

Let $\fa = {\rm span}_{{\mathbb R}}\{\i H_\alpha: \alpha \in \Sigma\}$, and let
$\ft = \i\fa$. We note that since
\[
\fh^{-\tau_v}=\fh^{-\gamma_d\circ\theta} = \ft^{-\gamma_d} + \fa^{\gamma_d},
\] 
the Lagrangian complement $\fl_d$ of $\fg_v$ depends only on $d$, and in the case when 
$d=1$, we have $\fl_d=\fa+\fn$. Note that $\fh^{\tau_v} = \fh^{\gamma_d\circ\theta}
= \ft^{\gamma_d} + \fa^{-\gamma_d}$ also depends only on $d$.

\begin{remark}
\label{rem_inner}
{\em
Recall  \cite[Definition 6.10]{a-v:l} that two real forms $\tau_1$ and $\tau_2$
are said to be in the same {\it inner class} if there exists $g \in {\rm Int}(\fg)$,  
the adjoint group 
of $\fg$, such that $\tau_1 = {\rm Ad}_g \tau_2$. Inner classes of real forms
are in one-to-one correspondence with involutive automorphisms of the Dynkin diagram of $\fg$  
\cite[Proposition 6.12]{a-v:l}. Let $d$ be an involutive automorphism of $D(\fg)$. 
Then as $v$ runs over the 
collection of all Vogan diagrams with $d$ as the diagram automorphism, the real form
$\fg_v$ runs over all ${\rm Int}(\fg)$-conjugacy classes of real forms of $\fg$ in the inner
class corresponding to $d$. 
}
\end{remark}

\section{The Poisson structure $\Pi_v$ on $X$.}
\label{sec_Piv}

Let $\fg$ be a complex semi-simple Lie algebra, and let $X$ be the variety of all Borel 
subalgebras of $\fg$. We keep the notation from Section \ref{sec_gv}. 
Let $v$ be a Vogan diagram for $\fg$ and $\fg_v= \fg^{\tau_v}$
be the real form of $\fg$ constructed in 
Section \ref{sec_gv}. Let $G$ be the connected and simple connected Lie group with
Lie algebra $\fg$. Without any risk of confusion, we shall also denote by $\tau_v$ the lift 
of $\tau_v$ from $\fg$ to $G$, and we set $G_v = G^{\tau_v}$. It follows 
from \cite[Theorem 8.2, p. 320]{Helga} that the group $G_v$ is connected.

In this section, we will start with a Vogan diagram $v$ for $\fg$ and
define a Poisson structure $\Pi_v$ on $X$ such that every $G_v$-orbit in $X$ is a Poisson 
submanifold. This Poisson structure comes from an identification of $X$ with
the $G$-orbit through $\ft+\fn$ inside the variety ${\mathcal L}$
of Lagrangian subalgebras of $\fg$, which was studied in \cite{EL}. 
We now recall the relevant details. 

Set $n = \dim_{{\mathbb C}} \fg$ and let ${\rm Gr}_\R(n,\fg)$ be
the Grassmannian of real $n$-dimensional subspaces of $\fg$. The set ${\mathcal L}$ of all
Lagrangian subalgebras of $\fg$ is naturally a real subvariety of 
${\rm Gr}_\R(n,\fg)$. The natural action of $G$ on ${\rm Gr}_\R(n, \fg)$ 
gives rise to a Lie algebra anti-homomorphism $\kappa$ from $\fg$ to the Lie algebra
of vector fields on ${\rm Gr}_\R(n,\fg)$, whose extension from $\wedge^2 \fg$ to 
the space of bi-vector fields on ${\rm Gr}_\R(n,\fg)$ will also be denoted by $\kappa$.
Given a Lagrangian splitting $\fg=\fl_1+\fl_2$, we define the element 
 $R_{{\mathfrak l}_1, {\mathfrak l}_2}\in\wedge^2\fg$ by: 
\begin{equation}
\label{eq_R}
\la R_{{\mathfrak l}_1, {\mathfrak l}_2}, \, (x_1+\xi_1) 
\wedge ( x_2+ \xi_2)\ra =\la \xi_2, x_1\ra -\la \xi_1, x_2\ra, \ \ \ \ 
x_1, x_2\in \fl_1, \ \xi_1, \xi_2\in \fl_2,
\end{equation}
where $\la \, , \, \ra = {\rm Im} \ll \, , \, \gg$. 
Set $\Pi_{{\mathfrak l}_1, {\mathfrak l}_2} = \frac{1}{2} \kappa(R_{{\mathfrak l}_1, 
{\mathfrak l}_2})$. 
Clearly, 
$\Pi_{{\mathfrak l}_1, {\mathfrak l}_2}$ is tangent to every $G$-orbit in 
${\rm Gr}_\R(n, \fg)$, so it is tangent to ${{\mathcal L}}$.

\begin{theorem}\cite[Theorems 2.14 and 2.18]{EL} The bi-vector field 
$\Pi_{{\mathfrak l}_1, {\mathfrak l}_2}$ restricts to a 
Poisson structure on ${\mathcal L}$. If $L_1$ and $L_2$ are the 
connected subgroups of $G$ with Lie algebras $\fl_1$ and $\fl_2$ respectively, then
all the $L_1$- as well as $L_2$-orbits in ${{\mathcal L}}$ are Poisson submanifolds
with respect to $\Pi_{{\mathfrak l}_1, {\mathfrak l}_2}$. 
\end{theorem}

For $\fl \in {{\mathcal L}}$, let $\fn(\fl)$ be the normalizer subalgebra
of $\fl$ in $\fl_1$. Let $\fm(\fl)$ be the annihilator of $\fn(\fl)$ in $\fl$, 
i.e. $\fm(\fl) = \{x \in \fl: \la x, y \ra = 0 \ \forall
y \in \fn(\fl)\} \subset \fl$, and let ${{\mathcal V}}(\fl) = \fn(\fl) + \fm(\fl)$.

\begin{proposition}\cite[Theorem 2.21]{EL} \cite[Corollary 7.3]{LuDuke}
\label{prop_rank}
For each $\fl \in {{\mathcal L}}$, the space ${{\mathcal V}}(\fl)$ is a Lagrangian
subalgebra of $\fg$. The co-dimension of the symplectic leaf of 
$\Pi_{{\mathfrak l}_1, {\mathfrak l}_2}$ through $\fl$ in the orbit $L_1 
\cdot \fl$ is equal to $\dim( {{\mathcal V}}(\fl) \cap \fl_2)$.
\end{proposition}

\begin{nota}
{\em Let $v$ be a Vogan diagram for $\fg$. 
We denote by $\Pi_v$ the Poisson structure on ${{\mathcal L}}$ defined by the 
Lagrangian splitting $\fg = \fg_v + \fl_d$ in Lemma \ref{lem_splitting}.
Let $H$, $N$, and $B$ be respectively the connected subgroups of $G$ with Lie algebras
$\fh$, $\fn$, and $\fb = \fh + \fn$, so $B = HN$. Identify the $G$-orbit through 
$\ft+ \fn \in {{\mathcal L}}$ with $G /B \cong X$. The induced Poisson structure 
on $X$ will also be denoted by $\Pi_v$.
Let 
$ H^{-\gamma_d\circ\theta} = \{h \in H: \gamma_d\circ\theta(h) = h^{-1}\}$ and let 
$L_d = H^{-\gamma_d\circ\theta}N$. By the Bruhat lemma, orbits of $L_d$ in 
$X \cong G/B$, which are the same as the $N$-orbits in $X$, are labeled by the elements 
in the Weyl group $W$ of $\Delta$. We refer to these $N$-orbits as the Bruhat cells in $X$. }
\end{nota}

By \cite[Theorem 2.18]{EL}, we have

\begin{proposition}
\label{prop_orbits-poisson-submanifolds}
Each $G_v$-orbit in $X$ as well as each Bruhat cell in $X$ is a Poisson submanifold 
with respect to $\Pi_v$.
\end{proposition}
 
When $v$ is the Vogan diagram with $d = 1$ and no vertex painted, we have 
$\tau_v = \theta$, so
$\fg_v = \fk$. The Poisson structure $\Pi_v$ in this case was first introduced in 
\cite{LW} and \cite{Soibel}, and it has the property that its symplectic leaves
are precisely the Bruhat cells (hence the name ``Bruhat Poisson structure" in
\cite{LW}). In \cite{EL} and \cite{LuCoord} this Poisson structure was related to some earlier work of
Kostant \cite{ko:63} and of Kostant-Kumar \cite{K-K:schubert} on the Schubert calculus on $X$.

The splitting $\fg = \fg_v + \fl_d$ naturally defines a Lie bialgebra structure on $\fg_v$
and therefore a Poisson Lie group structure on $G_v$ \cite{LW}. All the $G_v$-orbits 
in ${{\mathcal L}}$ become $G_v$-Poisson homogeneous spaces \cite{EL, LuDuke}. 
We remark that in \cite{a-j:real}, Andruskiewitsch and Jancsa classified  
non-triangular Lie bialgebra structures on $\fg_v$ using Belavin-Drinfeld triples. 
The one defined by the splitting $\fg = \fg_v + \fl_d$ comes from
the standard Belavin-Drinfeld triple. We refer to \cite{a-j:real} for details. 
 
\medskip \noindent {\bf Example.} Here we take 
$\fg={\mathfrak s}{\mathfrak l} (2, \C)$ and 
\[
\fg_v=\su(1,1) = \left\{ \left(\begin{array}{cc} \i x & y+\i z \\ y-\i z & -\i x 
\end{array} \right): \, x, y, z \in \R\right\}.
\]
Then $d = 1$ and $\fl_d = \fa + \fn$ consists of upper triangular matrices
in $\sl(2, \C)$ with real diagonal entries. Identify $G/B$ with $\P^1$ via the action 
\[
\left( \begin{array}{cc} a & b \\ c & d \end{array} \right) \cdot [w_0 : w_1] = 
[aw_0 + bw_1 : cw_0 + dw_1]
\]
of $G$ on $\P^1$
and by taking $[1 : 0] \in \P^1$ as the basepoint. There are two Bruhat cells: 
the zero-dimensional basepoint $[1 : 0]$, and the other being the rest:
\[
U_1 = \P^1 \backslash \{[1 : 0]\} = \{[w_0 : w_1], \ w_1 \neq 0\}.
\]
In terms of the holomorphic coordinate $z$ on $U_1$ given by $z=w_0/w_1$
the Poisson structure $\Pi_v$, up to a scalar multiple, is given by: 
\[
\Pi_v =  \i (1-|z|^2)\frac{\partial}{ \partial{z}} \wedge
\frac{\partial}{\partial{\bar z}}.
\]
Setting $u = 1/ z$, we see that in the $u$-coordinate on the open set
\[
U_0 = \{[w_0 : w_1] \in \P^1, w_0 \neq 0\} = \{[1 : u], u \in \C\},
\]
 we have
\[
\Pi_v=\i (|u|^2-1)|u|^2 \frac{\partial}{\partial{u}} \wedge
\frac{\partial}{\partial{\bar u}}.
\]
Thus $\Pi_v$ vanishes precisely at the basepoint $[1 : 0]$ and at every point of the form
$[z : 1]$ with $|z|=1$. If we identify $\P^1$ with the unit sphere $S^2$ in 
$\R^3$ via:
\begin{equation}
\P^1 \longrightarrow S^2: \, \, [w_0, w_1] \longmapsto 
\left( \frac{2{\rm Re}(w_0\overline{w_1})}{|w_0|^2+|w_1|^2}, \, \, 
\frac{2{\rm Im}(w_0\overline{w_1})}{|w_0|^2+|w_1|^2}, \, \, \frac{|w_0|^2 - 
|w_1|^2}{|w_0|^2+|w_1|^2}\right)\ ,
\end{equation}
then we see that $\Pi_v$ vanishes at the ``North pole" $(0, 0, 1)$ and at every point on the
Equator $x_3 = 0$. Under this identification, there are exactly three
orbits of ${\rm SU}(1, 1)$ on $S^2$: the Northern hemisphere, the Equator,
and the Southern hemisphere. Each one of these three orbits is clearly a Poisson
submanifold.
\medskip

\section{Symplectic leaves of $\Pi_v$ in $X$.} 

Suppose that ${{\mathcal O}}$ is a $G_v$-orbit in $X$ and ${{\mathcal C}}$ is 
a Bruhat cell such that ${{\mathcal O}}\cap {{\mathcal C}} \neq \emptyset$.
Since $\fg = \fg_v + \fl_d$, ${{\mathcal O}}$ and ${{\mathcal C}}$ intersect transversally.
By Proposition 
\ref{prop_orbits-poisson-submanifolds}, ${{\mathcal O}}\cap 
{{\mathcal C}}$ is a Poisson submanifold of $\Pi_v$. In this section we show that 
$({{\mathcal O}}\cap {{\mathcal C}}, \Pi_v)$ is a regular Poisson manifold and we compute
 the dimension of its symplectic leaves.

It is well-known \cite{wolf1} that there are only finitely many $G_v$-orbits in $X$. We
first recall from \cite[Section 6]{R-S:orbits} some facts about these orbits.  

Let $N_G(\fh)$ be the normalizer subgroup of $\fh$ in $G$. Set
\[
{{\mathcal Z}} = \{g \in G: \, g^{-1} \tau_v(g) \in N_G(\fh)\}.
\]
Then $H$ acts on ${{\mathcal Z}}$ from the right by right multiplication,
and $G_v$ acts on ${{\mathcal Z}}$ from the left by left multiplication. Let $Z$ be
the double coset space
\[
Z = G_v \backslash {{\mathcal Z}} / H.
\]
For each $z  \in Z$, choose any $g_z \in {{\mathcal Z}}$ in the double coset $z$ and 
define ${{\mathcal O}}_z$ to be the $G_v$-orbit in $X$ through $g_zB \in X \cong G/B$.
Clearly, ${{\mathcal O}}_z$ is independent of the choice of $g_z$. 
According to \cite[Theorem 6.1.4]{R-S:orbits}, the map 
$z \mapsto {{\mathcal O}}_z$ is a  one-to-one correspondence between the  set $Z$ and 
the set of $G_v$-orbits in $X$. Let $W = N_G(\fh)/H$ be the Weyl group. Thus we also have the map 
\[
\varphi: \, \, Z \longrightarrow W: \, \, z = G_v g_zH  \longmapsto g_{z}^{-1} \tau_v(g_z) H \in W.
\]
According to \cite[Theorem 6.4.2]{R-S:orbits}, the codimension of the 
$G_v$-orbit ${{\mathcal O}}_z$ in $X$ equals $l(\varphi(z))$, where $l$ is the length function on the Weyl group $W$.
We also introduce the map: 
\[
\sigma_z = \varphi(z) \tau_v: \, \, \fh \longrightarrow \fh.
\]
For any $g_z$ in the double coset $z$,  we also have
$\sigma_z = {\rm Ad}_{g_z}^{-1}\circ\tau_v\circ{\rm Ad}_{g_z}$, so $\sigma_z$ is an involution. 

Assume now that $z \in Z$ and $w \in W$ are such that ${{\mathcal O}}_z \cap {{\mathcal C}}_w 
\neq \emptyset$, where ${{\mathcal C}}_w$ is the Bruhat cell in $X$ corresponding to $w$, 
i.e. the $N$-orbit through $w \in G/B$. Then $\dim_{{\mathbb R}} {{\mathcal C}}_w = 2l(w)$, 
and since ${{\mathcal O}}_z$ and ${{\mathcal C}}_w$ intersect transversally, we have
\[
\dim({{\mathcal O}}_z \cap {{\mathcal C}}_w) = 2 l(w) - l(\varphi(z)).
\]
Define now 
\[
\delta_{z, w} = \dim (\fh^{w \sigma_z w^{-1}} \cap \fh^{-\tau_v}).
\]

\begin{theorem}
\label{thm_dim}
Each symplectic leaf in the intersection 
${\mathcal O}_z \cap {{\mathcal C}}_w$ has dimension equal to
$$
{\rm dim} ({{\mathcal O}}_z \cap {{\mathcal C}}_w) - \delta_{z,w} = 
2l(w) - l(\varphi(z)) - \delta_{z,w}.
$$
\end{theorem}

\proof We use Proposition \ref{prop_rank} to compute dimensions of
the symplectic leaves in ${{\mathcal O}}_z \cap {{\mathcal C}}_w$.
Let $x = g_zB \in X$ be a point in ${{\mathcal O}}_z \cap 
{{\mathcal C}}_w$, where $g_z \in {{\mathcal Z}}$ lies in the double coset
$z$.  
Let $\fl_x = {\rm Ad}_{g_z} (\ft + \fn) \in {{\mathcal L}}$. 
Let  $\fn(\fl_x) = \fg_v \cap {\rm Ad}_{g_z}(\fh + \fn)$ 
be the normalizer subalgebra of $\fl_x$ in $\fg_v$, $\fm(\fl_x)$ the annihilator
subspace of $\fn(\fl_x)$ in $\fl_x$, and ${{\mathcal V}}(\fl_x) = 
\fn(\fl_x) + \fm(\fl_x)$. We claim that ${{\mathcal V}}(\fl_x) = 
{\rm Ad}_{g_z}(\fh^{\sigma_z} + \fn)$. Indeed, it follows from the definition  of
$\sigma_z$ that $${\rm Ad}_{g_z}(\fh^{\sigma_z}) \subset \fg_v \cap 
{\rm Ad}_{g_z}(\fh + \fn) = \fn(\fl_x)\ .$$ It is also clear that 
${\rm Ad}_{g_z} \fn \subset \fm(\fl_x)$, so 
$${\rm Ad}_{g_z}(\fh^{\sigma_z} + \fn) \subset \fn(\fl_x) + \fm(\fl_x) = {{\mathcal V}}(\fl_x)\ .$$
Since both  ${\rm Ad}_{g_z}(\fh^{\sigma_z} + \fn)$ and  ${{\mathcal V}}(\fl_x)$ have the same 
dimension, they must coincide.

Let now $S_x$ be the symplectic leaf of $\Pi_v$ in $X$ through $x$. By Proposition \ref{prop_rank}, the
codimension of $S_x$ in ${{\mathcal O}}_z$ is equal to $\dim ({{\mathcal V}}(\fl_x) \cap
\fl_d)$. Let $\dot{w} \in N_G(\fh)$ be a representative of $w$ in $K$. Since $x \in 
{{\mathcal C}}_w$, there exist $n \in N$ and $b \in B$ 
such that $g_z = n \dot{w}b$. Then we have
\beqa
{{\mathcal V}}(\fl_x) \cap \fl_d  & = & \left({\rm Ad}_{n \dot{w}b} (\fh^{\sigma_z} + \fn)\right)
\cap (\fh^{-\tau_v} + \fn) \\
& = &  {\rm Ad}_{n}\left( ({\rm Ad}_{\dot{w}} (\fh^{\sigma_z} + \fn)) \cap 
(\fh^{-\tau_v} + \fn) \right)\\
& = & {\rm Ad}_{n}\left( \fh^{w \sigma_z w^{-1}} \cap \fh^{-\tau_v}
+ ({\rm Ad}_{\dot{w}} \fn) \cap \fn \right),
\eeqa
where in the last line we have the direct sum of vector spaces.
Since $$\dim ({\rm Ad}_{\dot{w}} \fn) \cap \fn  = \dim_\R X - \dim_\R {{\mathcal C}}_w\ , $$ we have
$$\dim ({{\mathcal V}}(\fl_x) \cap \fl_d) = \delta_{z,w} + \dim_\R X - \dim_\R {{\mathcal C}}_w\ ,$$ 
and thus
\[
\dim S_x = \dim {{\mathcal O}}_z -\dim ({{\mathcal V}}(\fl_x) \cap \fl_d) = \dim({{\mathcal O}}_z
\cap {{\mathcal C}}_w) - \delta_{z,w}.
\]
\qed  

Note, that the number $\delta_{z,w}$ depends only on $d$ and the two Weyl
group elements $\varphi(z)$ and $w$. Define $d: W \to W$ by $d(w) = \gamma_d w \gamma_d$. 
Following \cite{R-S:orbits}, we say that $w \in W$ is a {\it $d$-twisted involution} if
$d(w) = w^{-1}$. Denote by ${{\mathcal I}}_d$ the set of all
$d$-twisted involutions in $W$. Clearly, every $\varphi(z)$ is in ${{\mathcal I}}_d$. The 
Weyl group $W$ acts
on ${{\mathcal I}}_d$ by 
$$w_1 \ast w = w_{1} w d(w_{1}^{-1}) \ \ {\rm for} \ w_1 \in W, \ \ 
{\rm and}\ w \in {{\mathcal I}}_d \ ,$$
and the set $\varphi(Z) \subset {{\mathcal I}}_d$ is $W$-invariant.
In fact, the $W$-action on  $G/H$, given by $w \cdot gH = gw^{-1}H$, 
commutes with the left action of $G_v$ by left multiplication, and thus induces
a left action of $W$ on $Z$, which we denote by $w \cdot z$ for $w \in W$ and $z \in Z$.
It is also easy to see that $\varphi: Z \to W$ is $W$-equivariant, i.e.
$\varphi(w \cdot z) = w \ast \varphi(z)$ for all $w \in W$ and $z \in Z$. Similarly,
the involution $\tau_v: G \to G$ gives rise to an involution on $Z$ which depends only on
$d$. Denote this involution by $z \to d(z)$. Then we also have $\varphi(d(z)) = d \varphi(z) =
\varphi(z)^{-1}$.
As maps on $\fh$, we see that $w \sigma_z w^{-1} = (w \ast \varphi(z)) \tau_v$. Thus we also
have:
\[
\delta_{z, w} = {\rm dim} ( \fh^{(w \ast \varphi(z)) \tau_v} \cap \fh^{-\tau_v}).
\]

\begin{corollary} 1) When $w \ast \varphi(z) =1$, 
symplectic leaves of $\Pi_v$ in ${{\mathcal O}}_z \cap {{\mathcal C}}_w$ are
precisely its connected components.

2) Every open orbit ${{\mathcal O}}_z$ has an open symplectic leaf
${{\mathcal O}}_z \cap {{\mathcal C}}_{w_0}$, where $w_0$ is the longest element in $W$;

3) If $d = 1$, symplectic leaves in an open orbit ${{\mathcal O}}_z$ are 
precisely the connected components of intersections  of Bruhat cells with ${{\mathcal O}}_z$.  
\end{corollary}

\proof 1) When $w \ast \varphi(z) = 1$, we have $\delta_{z,w} = 0$, so every symplectic
leaf in ${{\mathcal O}}_z \cap {{\mathcal C}}_w$ is open in ${{\mathcal O}}_z \cap {{\mathcal C}}_w$.

2) Since ${\mathcal C}_{w_0}$ is dense in $X$, it intersects with every open orbit
${{\mathcal O}}_z$. Since an orbit ${{\mathcal O}}_z$  is open if and only if
$\varphi(z) = 1$, statement 2) follows from 1) and the fact that $w_0$ commutes with $d$.
The fact that ${\mathcal C}_{w_0}\cap {\mathcal O}_z$ is connected follows from the 
observation that ${\mathcal O}_z$ is a connected open complex submanifold of $X$ and thus 
${\mathcal O}_z\cap(X\setminus{\mathcal C}_{w_0})$ is a divisor in ${\mathcal O}_z$. 

3) follows directly from 1).
\qed

\bigskip

Consider now the group $H^{\tau_v} =H \cap G_v$. Since the centralizer of ${\fh}^{\tau_v}$
in $G_v$ also centralizes $\fh$, we see that $H^{\tau_v}$ is the Cartan subgroup of 
$G_v$ corresponding to the Cartan subalgebra ${\fh}^{\tau_v}$. Then according to 
\cite[Proposition 7.90]{Knapp} the group $H^{\tau_v}$ is connected.

The Poisson structure $\Pi_v$ on
$X$ is $H^{\tau_v}$-invariant. Indeed, let $R \in \wedge^2 \fg$ be the element given in
(\ref{eq_R}) for $\fl_1 = \fg_v$ and $\fl_2 = \fl_d$. We can also represent 
$R$ as $R = \sum_{i} \xi_i \wedge y_i$, where $\{y_i\}$ is a basis of $\fg_v$,
and $\{\xi_i\}$ is the dual basis of $\fl_d$ with respect to the pairing between
$\fg_v$ and $\fl_d$ given by $\la \, , \, \, \ra$, the imaginary part of the Killing form
on $\fg$. If $h \in H^{\tau_v}$, then $\{{\rm Ad}_h y_i\}$ is a basis of $\fg_v$, and
$\{{\rm Ad}_h \xi_i\}$ is its dual basis. Thus ${\rm Ad}_h R = R$. 
 
Assume now that $z \in Z$ and $w \in W$ are such that 
${{\mathcal O}}_z \cap {{\mathcal C}}_w \neq \emptyset$. Clearly, $H^{\tau_v}$ leaves 
${{\mathcal O}}_z \cap {{\mathcal C}}_w$ invariant. Since the Poisson structure $\Pi_v$ is 
$H^{\tau_v}$-invariant, if $S_x$ is the symplectic leaf of $\Pi_v$ through $x$, then
$hS_x:=\{hx_1: x_1 \in S_x\}$ is the symplectic leaf of $\Pi_v$ through $hx$.  
Define: $$F_x := \bigcup_{h \in H^{\tau_v}} hS_x\ .$$ 

\begin{proposition}
\label{prop_components}
For any $x \in X$, the set $F_x$ is a connected component of 
${{\mathcal O}}_z \cap {{\mathcal C}}_w$. 
\end{proposition}

\proof It is easy to see that if $F_{x_1} \cap F_{x_2} \neq
\emptyset$, then $F_{x_1} = F_{x_2}$. The statement would follow once we prove
that $F_x$ is an open subset of ${{\mathcal O}}_z \cap {{\mathcal C}}_w$ for each $x$.

Let $x = g_zB \in {{\mathcal O}}_z \cap {{\mathcal C}}_w$ with $g_z \in {{\mathcal Z}}$ in
the double coset $z$. For $y \in \fh^{\tau_v}$, let $X_y$ be the 
vector field on $X$ generating the action of $\exp(ty) \in H_{0}^{\tau_v}$ on $X$.
We claim that $X_y(x) \in T_{x}S_x$ if and only if $y \in p(\fh^{(w \ast \varphi(z)) \tau_v})$, where
$p: \fh \to \fh^{\tau_v}$ is the projection with respect to the decomposition 
$\fh = \fh^{\tau_v} + \fh^{-\tau_v}$. Assume the claim. Then since the kernel of the
map $p: \fh^{(w \ast \varphi(z)) \tau_v} \to \fh^{\tau_v}$ has dimension 
$\dim (\fh^{(w \ast \varphi(z)) \tau_v} \cap \fh^{-\tau_v}) = \delta_{z, w}$,
the image of the map
\[
J_x: \, \, \fh^{\tau_v} \longrightarrow T_x {{\mathcal O}}_z / T_xS_x: \, \, y \longmapsto 
X_y(x) + T_xS_x
\]
has dimension equal to $\dim (\fh^{\tau_v}) - \dim (\fh^{(w \ast \varphi(z)) \tau_v}) 
+ \delta_{z, w} = \delta_{z, w}$. Thus $J_x$ is onto, and the $H_{0}^{\tau_v}$-orbit
in ${{\mathcal O}}_z \cap {{\mathcal C}}_w$ through $x$ is transversal to the symplectic leaf
$S_x$. It follows that $F_x$ is open in ${{\mathcal O}}_z \cap {{\mathcal C}}_w$.

It remains to prove the claim. Denote also by $p: \fg \to \fg_{v}$
the projection with respect to
the decomposition $\fg = \fg_v + \fl_d$, and let $q$ be the projection
$q: \fg_v \to \fg_v / \fg_v \cap {\rm Ad}_{g_z} \fb \cong T_x {{\mathcal O}}_z$.
Then by \cite[Corollary 7.3]{LuDuke}, we have 
$T_xS_x = (q \circ p) ({{\mathcal V}}(\fl_x))$, where,
as in the proof of Theorem \ref{thm_dim}, 
${{\mathcal V}}(l_x)={\rm Ad}_{g_z} (\fh^{\sigma_z} + \fn)$.
Let $y \in \fh^{\tau_v}$. If $X_y(x) \in T_xS_x$, then there exist $y_1 \in \fl_d$
and $y_2 \in \fg_v$ with $y_1 + y_2 \in {{\mathcal V}}(\fl_x)$ such that 
$y-y_2 \in \fg_v \cap {\rm Ad}_{g_z} \fb \subset {{\mathcal V}}(\fl_x)$. Thus
$y+y_1 = y-y_2 + y_1 + y_2 \in {{\mathcal V}}(\fl_x)$. Write $y_1 = \xi_1 + u_1$, where
$\xi_1 \in \fh^{-\tau_v}$ and $u_1 \in \fn$. Then there exist $\xi_2 \in \fh^{\sigma_z}$ and
$u_2 \in \fn$ such that $y+\xi_1 + u_1 = {\rm Ad}_{g_z}(\xi_2 + u_2)$.
Write $g_z = n \dot{w} b$, where $n \in N, b \in B$, and $\dot{w}$ is a representative of
$w$ in $K$. Write ${\rm Ad}_{n^{-1}} (y+\xi_1+u_1) = y + \xi_1 + u_{1}^{\prime}$ and
${\rm Ad}_b (\xi_2 + u_2) = \xi_2 + u_{2}^{\prime}$, where $u_{1}^{\prime}, u_{2}^{\prime} 
\in \fn$. Then we have
\[
y+\xi_1 + u_{1}^{\prime} = {\rm Ad}_{\dot{w}} (\xi_2 + u_{2}^{\prime}).
\]
Since $y+\xi_1, {\rm Ad}_{\dot{w}} \xi_2 \in \fh$ and 
$u_{1}^{\prime}, {\rm Ad}_{\dot{w}} u_{2}^{\prime} \in \fn + \fn_-$, where $\fn_{-}
=\theta(\fn)$, we have $y + \xi_1 = {\rm Ad}_{\dot{w}} \xi_2 \in 
\fh^{(w \ast \varphi(z)) \tau_v}$. Thus $y \in p(\fh^{(w \ast \varphi(z)) \tau_v})$. 
Conversely, if $y \in \fh^{\tau_v}$ is such that $y + \xi_1 \in
\fh^{(w \ast \varphi(z)) \tau_v}= {\rm Ad}_{\dot{w}} \fh^{\sigma_z}$ 
for some $\xi_1 \in \fh^{-\tau_v}$, write $y + \xi_1 = {\rm Ad}_{\dot{w}} \xi_2$
for $\xi_2 \in \fh^{\sigma_z}$. Let ${\rm Ad}_{b^{-1}} \xi_2 =
\xi_2 + u_2$ for some $u_2 \in \fn$. We then have
\[
{\rm Ad}_n (y + \xi_1) = {\rm Ad}_{n \dot{w} b} (\xi_2 + u_2) \in {{\mathcal V}}(\fl_x).
\]
On the other hand, let ${\rm Ad}_n(y + \xi_1) = y + \xi_1 + u_1$ with $u_1 \in \fn$.
We see that $y= p ({\rm Ad}_n(y+\xi_1))$ so $X_y(x) \in T_xS_x$. 
\qed

\section{Invariant Poisson cohomology of open orbits.}

Let ${{\mathcal O}}_z$ be a $G_v$-orbit in $X$ equipped
with the Poisson structure $\Pi_v$. Then $({{\mathcal O}}_z, \Pi_v)$ 
is a Poisson homogeneous space for the Poisson Lie group $G_v$.
 The $G_v$-invariant Poisson 
cohomology of $({{\mathcal O}}_z, \Pi_v)$, denoted by $H^{\fd}_{\Pi_v, G_v}({{\mathcal O}}_z)$,
is defined as the cohomology of the cochain complex $(\chi^{\fd}_{G_v}({{\mathcal O}}_z), 
\partial_{\Pi_v})$, where $\chi^{\fd}({{\mathcal O}}_z)^{G_v}$ is the space of all 
$G_v$-invariant complex multi-vector fields on ${{\mathcal O}}_z$, $d_{\Pi_v}(V) = 
[\Pi_v, V],$
and $[\cdot, \cdot ]$ is the Schouten bracket of the multi-vector fields. 

\begin{proposition}
When ${{\mathcal O}}_z$ is an open $G_v$-orbit in $X$, the $G_v$-invariant
Poisson cohomology $H^{\fd}_{\Pi_v, G_v}({{\mathcal O}}_z)$ is isomorphic to the
de Rham cohomology of $X$.
\end{proposition}

\proof As in the proof of Theorem \ref{thm_dim},  let
$x = g_z B \in X$ be an arbitrary point in ${{\mathcal O}}_z$, where
$g_z \in {{\mathcal Z}}$ is in the coset $z$, and let 
${{\mathcal V}}(\fl_x) = 
{\rm Ad}_{g_z}(\fh^{\sigma_z} + \fn)$. Since  ${{\mathcal O}}_z$ is open, 
the stabilizer subalgebra of $\fg_v$ at $x$ is $\fg_v \cap {{\mathcal V}}(\fl_x)=
{\rm Ad}_{g_z} (\fh^{\sigma_z})$.
By \cite[Theorem 7.5]{LuDuke}, the $G_v$-invariant
Poisson cohomology $H^{\fd}_{\Pi_v, G_v}({{\mathcal O}}_z)$ is isomorphic to the
relative Lie algebra cohomology of the Lie algebra ${{\mathcal V}}(\fl_x) \otimes \C$ 
relative
to the subalgebra $({\rm Ad}_{g_z}(\fh^{\sigma_z}))\otimes \C$. 
Thus the $G_v$-invariant Poisson cohomology is
isomorphic to the $\fh$-invariant part of the Lie algebra cohomology of the
direct sum Lie algebra $\fn \oplus \fn$ with coefficients in $\C$, which by Kostant's
theorem \cite{ko:63}, is isomorphic to the de Rham cohomology of $X$.
\qed

\section{Remarks.}
We have constructed a Poisson structure $\Pi_v$ on $X$ for each Vogan diagram $v$ for
$\fg$ (which is not necessarily normalized). In particular, each Bruhat cell ${{\mathcal C}}_w$ 
in $X$ carries the Poisson structure $\Pi_v$. It would be interesting to study connections between
the Poisson structures for different $v$. Especially interesting are the properties of
$\Pi_v$ that depend only on the inner class $d$ of the real form $\fg_v$.  We also remark that 
the Poisson structure $\Pi_v$ is defined on the whole variety ${{\mathcal L}}$ of
Lagrangian subalgebras of $\fg$. We have only been looking
at the restriction of $\Pi_v$ to a particular $G$-orbit, namely the $G$-orbit through
the Lagrangian subalgebra $\ft+\fn$. There are many other interesting
$G$-orbits in ${{\mathcal L}}$, such as the $G$-orbit through a given real form of $\fg$.
It would be interesting to study the properties of the Poisson structure $\Pi_v$ on 
these orbits as well as on their closures with respect to both the classical topology and the 
Zariski topology.

\section*{Acknowledgments.}
Initial ideas of the paper came from discussions with Sam Evens.  
The first author was supported by NSF grant DMS-0072520. The second author
was supported by NSF grants DMS-0105195 and DMS-0072551 and by the 
HHY Physical Sciences Fund at the University of Hong Kong. Both authors
are grateful to IHES for hospitality and the second author also wishes to thank 
Rencontres Math\'ematiques de Glanon 2003, where the paper was completed.


\end{document}